\documentclass[12pt]{article}
\usepackage{mathrsfs}
\usepackage{amsfonts}
\usepackage{amssymb}
\usepackage{amsmath}
\usepackage{framed}
\usepackage{graphicx}
\usepackage{hyperref}
\usepackage{xcolor}
= 0pt \leftmargin     = 1.25in \topmargin =0pt \headheight     = 0pt \headsep
= 0pt \topskip =0pt
\def\sqr#1#2{{\vcenter{\vbox{\hrule height.#2pt
              \hbox{\vrule width.#2pt height#1pt \kern#1pt \vrule width.#2pt}
              \hrule height.#2pt}}}}
\def\dbR{{\mathop{\rm l\negthinspace R}}}

\def\3n{\negthinspace \negthinspace \negthinspace }
\def\2n{\negthinspace \negthinspace }
\def\1n{\negthinspace }
\def\ns{\noalign{\smallskip} }
\def\ns{\noalign{\medskip} }

\def\dbR{{\mathop{\rm l\negthinspace R}}}

\def\={\buildrel \triangle \over =}

\newcommand{\norm}[1]{\left\Vert#1\right\Vert}

\def\l{\lambda}

\def\f{\varphi}

\def\O{\Omega}

\def\no{\noindent}

\def\ms{\medskip}
\def\bs{\bigskip}
\def\q{\quad}
\def\qq{\qquad}

\def\ns{\noalign{\smallskip}}
\def\ds{\displaystyle}

\def\Im{{\mathop{\rm Im}\,}}

\def\|{\Big |}
\def\({\Big (}
\def\){\Big )}
\def\[{\Big[}
\def\]{\Big]}
\def\be{\begin{equation}}
\def\bel{\begin{equation}\label}
\def\ee{\end{equation}}
\def\bea{\begin{eqnarray}}
\def\eea{\end{eqnarray}}
\def\bt{\begin{theorem}}
\def\et{\end{theorem}}
\def\bc{\begin{corollary}}
\def\ec{\end{corollary}}
\def\bl{\begin{lemma}}
\def\el{\end{lemma}}
\def\bp{\begin{proposition}}
\def\ep{\end{proposition}}
\def\br{\begin{remark}}
\def\er{\end{remark}}
\def\ba{\begin{array}}
\def\ea{\end{array}}
\def\bd{\begin{definition}}
\def\ed{\end{definition}}
\newtheorem{lemma}{Lemma}[section]
\newtheorem{remark}{Remark}[section]

\newtheorem{theorem}{Theorem}[section]
\newtheorem{corollary}{Corollary}[section]

\newtheorem{definition}{Definition}[section]

\newtheorem{proposition}{Proposition}[section]
\makeatletter
   
   \@addtoreset{equation}{section}
\makeatother
\begin{document}
 \title{\bf Inverse problems for the fourth order  Schr\"odinger equation on a finite domain
 \thanks{ This work was partially supported by
  the NSF of China under grant 11001018 and SRFDP (No.201000032006). This work has been completed while the author visited BCAM - Basque Center for Applied Mathematics and he acknowledges the hospitality and support of the Institute.
\ms} }
\author{ Chuang Zheng\thanks{ School of Mathematics,
 Beijing Normal University,
 100875 Beijing,
 China.    chuang.zheng@bnu.edu.cn.\ms 
}$\;$$\;$
  }
 \maketitle

\begin{abstract}
In this paper we establish a global Carleman estimate for the fourth order Schr\"odinger equation posed on a  $1-d$ finite domain.
The Carleman estimate is used to prove the Lipschitz stability for an inverse problem consisting in retrieving a stationary potential in the Schr\"odinger equation from boundary measurements.
\end{abstract}

\bs

\no{\bf AMS Subject Classifications}.
93B05, 
35Q40. 
\bs

\no{\bf Key Words}. Inverse problem, Fourth order Schr\"odinger equation,  Carleman estimate.

\section{Introduction}

 The fourth order Schr\"odinger equation arises in many scientific fields such as quantum mechanics, nonlinear optics and plasma physics, and has been intensively studied with fruitful references. The well-posedness and existence of the solutions has been shown (for instance, see \cite{MR2230468, MR2353631, MR2502523}) by means of the energy method and harmonic analysis.  In this paper, we are interested in the inverse problem for the fourth order Schr\"odingier equation posed on a finite interval.

To be more precise,  we consider the following fourth order Schr\"odinger equation on $\Omega=(0,1)$:
  \begin{equation}\label{maineqn}
  \left\{\begin{array}{ll}
  i u_t+u_{xxxx}+pu=0, &(t,x)\in (0,T)\times\O\\
  u(t,0)=u(t,1)=0, u_x(t,0)=u_x(t,1)=0, & t\in (0, T)\\
u(0,x)=u_0(x), & x\in \O.
   \end{array}\right.\end{equation}
For any initial data $u_0\in H^3(\Omega)\cap H_0^2(\O)$ and $p\in L^2(\O)$,  there exists a unique solution of \eqref{maineqn} $u\in C^1([0,T]; L^2(\O))\cap C([0,T]; H^3(\Omega)\cap H_0^2(\O))$(see, for instance, \cite{MR953547}).

The purpose of this paper is to determine  the potential $p=p(x), x\in \Omega$ by means of the boundary measurements. The problem we are interested can be stated as follows: is it possible to estimate $\norm{q-p}_{L^2(\Omega)}$, or better, a stronger norm of $q-p$, by a suitable norm of the derivatives of $u(q)-u(p)$ at the end point $x=1$ (or, at $x=0$) during the time interval $(0,T)$?

Recently, the inverse problem of the Schr\"odinger equations have been intensely studied (see \cite{MR1955903, MR2557958, MR2005314, MR2864510, MR1862200, MR2384776,  MR2672251} and the references therein). One of the main techniques is the Carleman estimate (\cite{MR1955903,  MR2864510,  MR2040524, MR2384776,  MR2672251}), which is also a powerful tool for the controllability and observability problems of PDEs.

However, for the higher order equations,  due to the increased complexity, there are few papers investigating the stability of the inverse problems via Carleman estimates.   In  \cite{MR1852002},  Zhang solves  the exact controllability of semilinear plate equations via a Carleman estimate of the second order Schr\"odinger operator.
Zhou (\cite{MR300183})  considers the observability results of the fourth order parabolic equation and Fu (\cite{MR2886515}) derives  the sharp observability inequality for the plate equation. In both papers, they show the Carleman estimates for the corresponding fourth order operators for $1-d$ cases, respectively.

To our knowledge, the result  of determination of a time-independent potential for the fourth order Schr\"odinger equation from the boundary measurements on the endpoint is new.  Furthermore, our work in this paper is the first one dealing with the Carleman estimate of the fourth order Schr\"odinger equation.

To begin with, we introduce a suitable weight function:
\bel{x_0}
\psi(x)=(x-x_0)^2, \qq x_0<0.
\ee
Let $\lambda\gg 1$ be a sufficiently large positive constant depending on $\Omega$. For $t\in (0,T)$ and following \cite{MR1750109}, we introduce the functions
\bel{ellvarphi}
\theta=e^l,\qq\f(t,x)=\ds \frac{e^{3\mu\psi(x)}}{t(T-t)}\qq\hbox{and}\qq l(t,x)= \l \frac{e^{3\mu\psi(x)}-e^{5\mu\norm{\psi}_\infty}}{t(T-t)}
\ee
with a positive constant $\mu$.  Denote by
 $$
 Pu=iu_t+u_{xxxx},\qq Q=(0,T)\times\Omega\q\hbox{and}\q \int_Q(\cdot)dxdt =\int_0^T\int_\Omega(\cdot)dxdt.
 $$
 We also introduce the set
\begin{equation*}
\begin{array}{ll}
&\ds        \mathcal{Z}\=    \big\{         Pu\in L^2(Q),  \;u\in L^2(0,T; H^3(\Omega)\cap H_0^2(\O)),\\
\ns
&   \ds        \qq\qq       u_{xx}(\cdot, 1)\in L^2(0,T),\; \; u_{xxx}(\cdot, 1)\in L^2(0,T)\big\}.
\ea
\end{equation*}
 The first main result is the following global Carleman estimate for system \eqref{maineqn}.
  %
 %
 \bt\label{Carl}
There exist three constants $\mu_0>1$, $C_0>0$ and $C>0$ such that for all $\mu\geq\mu_0$ and for all $\l\geq C_0 (T+T^2)  $,
\bel{Carest}
\ba{ll}
     &\ds \int_Q \(\l^7\mu^8\f^7\theta^2|u|^2+\l^5 \mu^6\f^5\theta^2|u_{x}|^2 +\l^3\mu^4\int_Q \f^3\theta^2|u_{xx}|^2+\l\mu^2\f \theta^2|u_{xxx}|^2\)dxdt\\
\ns\ns\ns
&\leq \ds C      \left(   \int_Q|\theta P u|^2dxdt
                         +\l^3\mu^3\int_0^T(\f^3\theta^2|u_{xx}|^2)(t,1) dt
                         +\l\mu \int_0^T (\f \theta^2|u_{xxx}|^2)(t,1) dt\right)
\ea
\ee
holds true for all $u\in\mathcal{Z}$, where the constants $\mu_0, C_0$ and $C$ only depend on $x_0$.
 \et
\begin{remark}
Note that for simplicity, we give the exact form of the function $\psi(x)$ in \eqref{x_0}. In fact,  the statement holds true for any function satisfying
$$
\psi\in C^4(\bar\Omega),\q\psi>0,\q \psi_x\neq 0\q\hbox{in}\q\bar\Omega\q\hbox{with}\q\bar\Omega\=[0,1], \qq\psi_x(0)>0, \psi_x(1)>0.
$$
It is worthy to mention that, by taking $x_0>1$,  one could switch the observation data in \eqref{Carest} to the left end-point $x=0$.  
\end{remark}
\begin{remark}\label{rem1.2}
\cite{MR2928966}  shows an observability inequality which estimates initial data  by the measurement of $\Delta u$  for a Schr\"odinger equation without the potential $q$  on $\Gamma_0=\{x\in \partial \Omega; (x-x_0)\cdot \nu(x)\geq 0\}$ using a multiplier identity and Holmgren's uniqueness  theorem. Observability inequalities are technically related to our inverse problem (see \cite{MR1671221}).  However,  the approach in \cite{MR2928966} can not be applied to our problem, even though there are less observability data are considered.

\end{remark}
\br
Note that the Carleman estimate \eqref{Carest} also can be applied to the controllability problems. In fact,  one can derive the exact controllability of the controlled fourth order semi-linear Schr\"odinger equations, with controls are given at the boundary point $x=1$.

\er

In  what  follows, we shall denote by $u^p$ the solution of the system \eqref{maineqn} associated with the potential $p$.

Following the standard procedure from the Carleman estimate to the inverse problem (see, for instance, \cite{MR2384776}), we answer the previous question with the following Theorem:
\begin{theorem}\label{invP}
Suppose that $p \in L^\infty(\Omega, \dbR)$, $ u_0\in L^\infty(\Omega)$ and $r>0$ are such that
\begin{itemize}
\item $u_0(x)\in\dbR$ or $iu_0(x)\in\dbR \;a.e.$ in $\Omega$,
\item $|u_0(x)|\geq r>0\; a.e. $ in $\Omega$, and
\item $u^p\in H^1(0,T; L^\infty(\Omega))$.
\end{itemize}
Then, for any $m\geq 0$, there exists a constant $C=C(m, \norm{u^p}_{H^1(0,T;L^\infty(\Omega))}, r)>0$ such that for any $q\in  L^\infty (\Omega, \dbR)$ satisfying
\begin{eqnarray}\label{12.13-eq1}
u_{xx}^p(t,1)-u_{xx}^q(t,1)\in
H^1(0,T)\qq\hbox{and}\qq
u_{xxx}^p(t,1)-u_{xxx}^q(t,1)\in H^1(0,T),
\end{eqnarray}
we have that
\begin{equation}\label{12.13-eq2}
\norm{p-q}^2_{L^2(\Omega)}\leq C
\(\norm{u_{xx}^p(\cdot,1)-u_{xx}^q(\cdot,1)}_{H^1(0,T)}^2+\norm{u_{xxx}^p(\cdot,1)-u_{xxx}^q(\cdot,1)}_{H^1(0,T)}^2\).
\end{equation}
\end{theorem}

\begin{remark}
By the classical regularity results for fourth
order Schr\"{o}dinger equations (see
\cite[Chapter 2]{MR1721989} for example), we
know that the $q$ which fulfills
\eqref{12.13-eq1} and \eqref{12.13-eq2} does
exist.
\end{remark}

The rest of the paper is organized as follows. In Section \ref{sec2}, we state a weighted point wise inequality for the fourth order Schr\"odinger  operator.  In Section \ref{sec3}, we establish a global Carleman estimate for a fourth order Schr\"odinger equation with a potential.  The proof of Theorem \ref{invP} is given in Section \ref{sec4}. Finally we list several comments and some open problems for the future work.

\section{A weighted point-wise estimate for the fourth order operator}\label{sec2}
In this section, we shall establish a weighted identity for $1$-d Schr\"odinger operator, which will pay an important role in the proof of the Carleman estimate \eqref{Carest}.
 %
 %
 \bt\label{Carleman}
Let $\Psi\in C^2(\dbR)$ and $v=\theta u$.  Then
\bel{PWI}
\ba{lll}
&&\ds|\theta P u|^2-A_x-B_t-\tilde a_0\theta (Pu\bar v+\overline{Pu} v)-6l_{xx}\theta(Pu\bar v_{xx}+\overline{Pu}v_{xx})\\
\ns\ns
&& \3n\3n=\ds |I_1|^2+|I_2|^2+ D(v\bar v_x-\bar v v_x)+6il_tl_{xx}(v\bar v_{xx} -\bar vv_{xx} )  +4il_{tx}(v_{xx}\bar v_x-\bar v_{xx}v_x)\\
\ns\ns
&&\ds +16l_{xx}|v_{xxx}|^2+(24l_x^2l_{xx}-24l_xl_{xxx}+48l_{xx}^2-20l_{xxxx})|v_{xx}|^2 \\
\ns\ns
&&\ds
\left\{\ba{rr}
\ds-(a_1a_2)_x-2(a_0-\Psi)a_2-4C_{41,xx}+3C_{24,x}\\
\ns\ns
\ds-a_{1,xxx}-\frac{3}{2}(a_{3,x}a_1)_x-3a_{3,x}a_0+2a_2\tilde a_0
\ea\right\}|v_{x}|^2 \\
\ns\ns
&&
\left\{
\begin{array}{rr}
+2(a_0-\Psi)\Psi-(a_1\Psi)_x-2a_0\tilde a_0-C_{24,xxx}+[(a_0-\Psi)a_2]_{xx}
\\
\ns
\ds-(a_1\tilde a_0)_x+(\frac{3}{2}a_{3,x} a_0-a_2\tilde a_0)_{xx}+C_{41,xxxx}+l_{tt}
\end{array}
\right\}
|v|^2,
\ea
\ee

where
\bel{a1234}
\left\{
\ba{ll}
 a_0\=\ds l_x^4-6l_x^2l_{xx}+3l_{xx}^2+4l_xl_{xxx}-l_{xxxx},\q  &a_2 \= 6 (l_x^2-l_{xx}),\\
    \ns
    \ns
   a_1\=-4(l_x^3-3l_xl_{xx}+l_{xxx}), \qq
  &a_3\=-4l_x,\\
  \ns\ns
  \tilde a_0=l_x^4-\Psi-2l_xl_{xxx}-3l_{xx}^2,\\
  \ns\ns
  C_{24}\=4l_x(l_x^4-2\Psi-2l_xl_{xxx}-3l_{xx}^2),\\
  \ns\ns\ds
  C_{41}\=-6l_x^2l_{xx}+6l_xl_{xxx}+6l_{xx}^2-l_{xxxx},
\ea
\right.
\ee
and
\begin{equation}\label{def_i}
I_1= i v_t+\Psi v+a_2 v_{xx}+ v_{xxxx}, \qq I_2=-i l_t v+ (a_0-\Psi) v+a_1v_{x}+a_3 v_{xxx},
\end{equation}
where $\Psi$ is a real value function in $C^2(\dbR)$. Moreover, we have
\bel{TermA}
\ba{lll}
&&\ds\3n\3n \3n \3n A= ia_3(v_{t}\bar v_{xx}-\bar v_{t} v_{xx})-\frac{i}{2}a_3(v_{xt}\bar v_{x}-\bar v_{xt} v_{x})+\frac{i}{2}a_{3,x}(v_{t}\bar v_{x}-\bar v_{t} v_{x})\\
\ns&&\ds+\frac{3}{2}a_{3,x}(v_{xxx}\bar v_{xx}+\bar v_{xxx}v_{xx})+a_1(v_{xxx}\bar v_{x}+\bar v_{xxx}v_{x})+C_{41}(v_{xxx}\bar v+\bar v_{xxx}v)\\
\ns
& &\ds+il_t(v_{xxx}\bar v-\bar v_{xxx}v)-il_t(v_{xx}\bar v_{x}-\bar v_{xx}v_{x})+(C_{24}-C_{41,x})(v_{xx}\bar v+\bar v_{xx}v)\\
\ns
& &\ds-il_{tx}(v_{xx}\bar v-\bar v_{xx}v)+i(l_{txx}+a_2l_t)(v_x\bar v-\bar v_x v)+\frac{i}{4}(2a_1-a_{3,xx}) (v_{ t}\bar v -\bar v_{  t} v )\\
\ns
& &\ds+[(a_0-\Psi)a_2-C_{24,x}-a_{1,x}-C_{41,xx}+\frac{3}{2}a_{3,x}a_0-a_2\tilde a_0](v_{x}\bar v +\bar v_{ x}v )\\
\ns
& &\ds+a_3|v_{xxx}|^2+(a_2a_3-\frac{3}{2}a_{3,x}a_3-\frac{3}{2}a_{3,xx}-a_1)|v_{xx}|^2\\
\ns
& &\ds+(a_1a_2+a_{1,xx}-C_{24}-2C_{41,x}+\frac{3}{2}a_{3,x}a_1)|v_{x}|^2\\
\ns
& &\ds+\left\{a_1\Psi+C_{24,xx}-C_{41,xxx}+[(a_0-\Psi)a_2]_x+a_1\tilde a_0-(\frac{3}{2}a_{3,x}a_0-a_2\tilde a_0)_{xx}\right\}|v|^2,\\
\ns\ns
&&\ds \3n\3n \3n \3n B=-l_t|v|^2 -\frac{i}{2}a_3(v_{x }\bar v_{xx}-\bar v_{x} v_{xx})  +\frac{i}{4}(2a_1-a_{3,xx}) (v \bar v_x -\bar v  v_x ),\\
\ns\ns
&&\ds\3n\3n \3n \3n D=2i(6l_x^2l_{xt}+6l_xl_{xx}l_t-6l_{xx}l_{xt}-3l_xl_{xxt}-3l_tl_{xxx}+l_{xxxt}).
\ea
\ee
 \et
%
%
{\bf Proof.}  We may assume that $u$ is sufficiently smooth.  Since $v=\theta u$ and notice the definitions of $a_i, i=0,1,2,3$ in \eqref{a1234}, it is esay to get
\begin{equation}\label{thetapu}
\theta P u=iv_t-il_tv+ v_{xxxx}+ a_0 v+a_1v_{x}+a_2 v_{xx}+a_3 v_{xxx}.
\ee
We divide $Pu$ into $I_1$ and $I_2$ as in \eqref{def_i}.
Multiplying $\theta P u$ by its conjugate we have
\bel{pu1}
|\theta P u|^2=|I_1|^2+|I_2|^2+(I_1\bar I_2+\bar I_1 I_2)\=|I_1|^2+|I_2|^2+\sum_{i,j=1}^4 I_{ij},
\ee
where $I_{ij}$ denotes the sum of the $i$-th term of $I_1$ times the $j$-th term of $\bar I_2$ in $I_1\bar I_2$ and its conjugate part in $\bar I_1 I_2$.

The computations will be treated in the following two parts.

{\bf Part I:} We compute $I_{1j}, j=1,2, 3, 4$. We first have
\bel{11}
I_{11}=-l_t(v_t\bar v+\bar v_t v)=-(l_t|v|^2)_t+l_{tt}|v|^2.
\ee

On the other hand, it is easy to verify that
\bel{13}\ba{rl}
I_{13}\=&ia_1(v_t\bar v_x-\bar v_tv_x)\\
\ns\ns
=&\ds\3n-\frac{i}{2}a_{1,x}(v_t\bar v-\bar v_tv)+\frac{i}{2}\{[a_1(v\bar v_x-\bar vv_x)]_t+[a_1(v_t\bar v-\bar v_tv)]_x-a_{1,t}(v\bar v_x-\bar vv_x)\}.
\ea\ee

Moreover,
\bel{14}\ba{rl}
I_{14}\=&ia_3(v_t\bar v_{xxx}-\bar v_tv_{xxx})\\
\ns\ns
=&\ds\3n-\frac{3i}{2}a_{3,x}(v_{t}\bar v_{xx}-\bar v_{t}v_{xx})-\frac{i}{2}a_{3,xx}(v_{t}\bar v_{x}-\bar v_{t}v_{x})\\
\ns\ns
&\ds-\frac{i}{2}[a_3(v_{x}\bar v_{xx}-\bar v_{x}v_{xx})]_t+\frac{i}{2}a_{3,t}(v_{x}\bar v_{xx}-\bar v_{x}v_{xx})\\
\ns\ns&\ds+\{ia_3(v_t\bar v_{xx}-\bar v_t v_{xx})-\frac{i}{2}a_3(v_{xt}\bar v_{x}-\bar v_{xt}v_{x})+\frac{i}{2}a_{3,x}(v_{t}\bar v_{x}-\bar v_{t}v_{x})\}_x.\\
\ea\ee

 By replacing $a_1$ in \eqref{13} by $a_{3,xx}$,  substituting it into the last term of \eqref{14}, we have
 \bel{141}\ba{rl}
I_{14}=&\ds\3n-\frac{3i}{2}a_{3,x}(v_{t}\bar v_{xx}-\bar v_{t}v_{xx})+\frac{i}{4}a_{3,xxx}(v_{t}\bar v-\bar v_{t}v)\\
\ns\ns
&\ds-\frac{i}{4}\{[a_{3,xx}(v\bar v_x-\bar vv_x)]_t+[a_{3,xx}(v_t\bar v-\bar v_tv)]_x-a_{3,xxt}(v\bar v_x-\bar vv_x)\}\\
\ns\ns
&\ds-\frac{i}{2}[a_3(v_{x}\bar v_{xx}-\bar v_{x}v_{xx})]_t+\frac{i}{2}a_{3,t}(v_{x}\bar v_{xx}-\bar v_{x}v_{xx})\\
\ns\ns&\ds+\{ia_3(v_t\bar v_{xx}-\bar v_t v_{xx})-\frac{i}{2}a_3(v_{xt}\bar v_{x}-\bar v_{xt}v_{x})+\frac{i}{2}a_{3,x}(v_{t}\bar v_{x}-\bar v_{t}v_{x})\}_x.\\
\ea\ee

 Set
 $$
 \tilde a_0\=\ds a_0-\Psi-\frac{1}{2}a_{1,x}-\frac{1}{4}a_{3,xxx}.
 $$
 Obviously, it is the coefficient of the term $i(v_{t}\bar v-\bar v_{t}v)$ in $\ds\sum_{j=1}^4 I_{1j}$. Taking the exact form of $a_0, a_1, a_3$ in \eqref{a1234} into account,   one can verifty that $\tilde a_0$ is exactly the one in \eqref{a1234}. Furthermore,
 \bel{15}\ba{rl}
\tilde a_0i(v_{t}\bar v-\bar v_{t}v)=&\ds\3n \tilde a_0 \theta (P u\bar v+\overline{Pu} v)-2a_0\tilde a_0|v|^2-a_1\tilde a_0(v_x\bar v+\bar v_x v)\\
\ns\ns
&\ds\3n-a_2\tilde a_0(v_{xx}\bar v+\bar v_{xx} v)-a_3\tilde a_0(v_{xxx}\bar v+\bar v_{xxx} v)-\tilde a_0(v_{xxxx}\bar v+\bar v_{xxxx} v).
\ea\ee
Meanwhile,  for the first term of $I_{14}$, recalling that $a_3=-4l_x$, we have
 \bel{16}\ba{rl}
&\ds\3n-\frac{3i}{2}a_{3,x}(v_{t}\bar v_{xx}-\bar v_{t}v_{xx})\\
\ns\ns
=&\ds\3n 6 l_{xx} \theta (P u\bar v_{xx}+\overline{Pu} v_{xx})+6il_{xx}l_t(v\bar v_{xx}-\bar v_{xx})\\
\ns\ns
&\ds\3n-6l_{xx}a_0(v_{xx}\bar v+\bar v_{xx} v)-6l_{xx}a_1(v_{xx}\bar v_x+\bar v_{xx} v_x)-12l_{xx}a_2|v_{xx}|^2\\
\ns\ns
&\ds\3n-6l_{xx}a_3(v_{xxx}\bar v_{xx}+\bar v_{xxx} v_{xx})-6l_{xx}(v_{xxxx}\bar v_{xx}+\bar v_{xxxx} v_{xx}).\\
\ea\ee

Summing up $I_{12}=i(a_0-\Psi)(v_{t}\bar v-\bar v_{t}v)$,  $I_{13}$ as \eqref{13} and $I_{14}$ as \eqref{141}, taking \eqref{15} and \eqref{16} into accout, we arrive at
\bel{I234}\ba{rl}
\ds\sum_{j=2,3,4}I_{1j}=&\{\cdot\}_x+\{\cdot\}_t+\tilde a_0\theta (Pu\bar v+\overline{Pu} v)+6l_{xx}\theta(Pu\bar v_{xx}+\overline{Pu}v_{xx})\\
%
&\ds+\frac{3}{2}a_{3,x}(v_{xxxx}\bar v_{xx}+\bar v_{xxxx} v_{xx}) -\tilde a_0(v_{xxxx}\bar v +\bar v_{xxxx} v)\\
\ns
& \ds+\frac{3}{2}a_{3,x}a_3(v_{xxx}\bar v_{xx}+\bar v_{xxx} v_{xx}) -a_3\tilde a_0(v_{xxx}\bar v +\bar v_{xxx} v) \\
\ns
& \ds+\frac{3}{2}a_{3,x}a_1(v_{xx}\bar v_{x}+\bar v_{xx} v_{x})  +(\frac{3}{2}a_{3,x}a_0-a_2\tilde a_0) (v_{xx}\bar v +\bar v_{xx} v)  \\
\ns
&  \ds+\frac{i}{2}a_{3,t} (v_{ x}\bar v_{xx}-\bar v_{x} v_{xx})+\frac{3}{2}i l_ta_{3,x} (v_{xx}\bar v-\bar v_{xx} v ) +3 a_{3,x }a_2|v_{xx}|^2  \\
\ns
& \ds +\frac{i}{4}(a_{3,xxt}-2a_{1,t}) (v \bar v_x-\bar v v_x )  + a_1\tilde a_0  (v_{x}\bar v +\bar v_{x} v) -2a_0\tilde a_0|v |^2,
\ea\ee
with
$$
\{\cdot\}_x =\left(\ba{ll}&\ds ia_3(v_{t}\bar v_{xx}-\bar v_{t} v_{xx})-\frac{i}{2}a_3(v_{xt}\bar v_{x}-\bar v_{xt} v_{x})\\
 \ns&\ds+\frac{i}{2}a_{3,x}(v_{t}\bar v_{x}-\bar v_{t} v_{x})+\frac{i}{4}(2a_1-a_{3,xx}) (v_{ t}\bar v -\bar v_{  t} v )\ea\right)_x
$$
and
$$
\{\cdot\}_t =\left( -\frac{i}{2}a_3(v_{x }\bar v_{xx}-\bar v_{x} v_{xx})  +\frac{i}{4}(2a_1-a_{3,xx}) (v \bar v_x -\bar v  v_x )\right)_t.
$$

\vskip 5mm

{\bf Part II:} We compute the rest of $I_{ij}$,  with some extra terms  coming from \eqref{I234}.

Set $C_{24}=a_3\Psi-a_3\tilde a_0,$ which is the same notation as in \eqref{a1234}. We have the following identity:
\bel{24}\ba{rl}
&\ds I_{24}-a_3\tilde a_0(v_{xxx}\bar v +\bar v_{xxx} v ) =C_{24}(v_{xxx}\bar v +\bar v_{xxx} v ) \\
\ns\ns
=&\ds
\left(\ba{rr}C_{24}(v_{xx}\bar v+\bar v_{xx} v)-C_{24}|v_x|^2\\
\ns-C_{24,x}(v_{x}\bar v+\bar v_{x} v)+C_{24,xx}|v |^2\ea\right)_x+3C_{24,x}|v_x|^2-C_{24,xxx}|v |^2.
\ea\ee
Consequently, it holds
\bel{i2j}\ba{rl}
&\ds\sum_{j=1}^4I_{2j} -a_3\tilde a_0(v_{x xx}\bar v +\bar v_{xxx } v ) \\
\ns
=&\ds 0+2(a_0-\Psi)\Psi v \bar v+a_1\Psi(v_{x }\bar v +\bar v_{x } v ) +C_{24}(v_{x xx}\bar v +\bar v_{xxx } v ) \\
\ns
=&\ds \{C_{24}(v_{xx}\bar v+\bar v_{xx} v)-C_{24,x}(v_{x}\bar v+\bar v_{x} v)-C_{24}|v_x|^2+(a_1\Psi+C_{24,xx})|v |^2\}_x\\
\ns\ns&\ds\qq\qq\qq\qq
+3C_{24,x}|v_x|^2+\{2(a_0-\Psi)\Psi-(a_1\Psi)_x-C_{24,xxx}\}|v |^2.
\ea\ee

Now we compute $I_{3j}, j=1,2,3,4.$  It holds
\bel{31}
I_{31} =ia_2l_t(v_{xx}\bar v -\bar v_{xx} v ) =\{ia_2l_t(v_{x}\bar v-\bar v_{x} v)\}_x-(ia_2l_{t})_x(v_{x}\bar v-\bar v_{x} v),
\ee
and
\bel{i3j}\ba{rl}
\3n\3n\ds\sum_{j=2,3,4}I_{3j}=&\3n(a_0-\Psi)a_2(v_{xx}\bar v+\bar v_{xx}v)+a_1a_2(v_{xx}\bar v_x+\bar v_{xx} v_x)+a_2a_3(v_{xx}\bar v_{xxx}+\bar v_{xx} v_{xxx})\\
\ns\ns=&\ds C_x + [(a_0-\Psi)a_2]_{xx}|v |^2-[(a_1a_2)_x+2(a_0-\Psi)a_2]|v_x|^2-(a_2a_3)_{x}|v_{xx}|^2,
\ea\ee
with
\be
C =  (a_0-\Psi)a_2(v_{x }\bar v+\bar v_{x} v)  +[(a_0-\Psi)a_2]_x|v |^2+a_1a_2|v_x|^2+ a_2a_3 |v_{xx}|^2.
\ee

For the term $I_{41}$, it holds:
\bel{I41}\ba{rcl}
I_{41} &=&il_t(v_{xxxx}\bar v -\bar v_{xxxx} v ) \\
\ns
&=&   [il_t(v_{xxx}\bar v-\bar v_{xxx} v-v_{xx}\bar v_x+\bar v_{xx} v_x)-il_{tx}(v_{xx}\bar v -\bar v_{xx} v )\\
\ns
&&+il_{txx}(v_{x}\bar v-\bar v_{x} v)]_x+2il_{tx}(v_{xx}\bar v_x-\bar v_{xx} v_x)-il_{txxx}(v_{x}\bar v-\bar v_{x} v).
\ea\ee
$I_{42}$ is considered with an extra term from $I_{14}$ as follows:
\bel{I42}\ba{rl}
\ds I_{42}-\tilde{ a}_0(v_{xxxx}\bar v +\bar v_{xxxx} v )\=&C_{41}(v_{xxxx}\bar v +\bar v_{xxxx} v )\\
\ns
=& E_x+2C_{41}|v_{xx}|^2-4C_{41,xx}|v_{x}|^2+C_{41,xxxx}|v|^2,
\ea\ee
 with
 \bel{E}\ba{rr}
  E=&\ds C_{41}(v_{xxx}\bar v+\bar v_{xxx} v)-C_{41}(v_{xx}\bar v_{x}+\bar v_{xx} v_{x})-C_{41,x}(v_{xx}\bar v +\bar v_{xx} v )\\
  \ns
  &\ds+C_{41,xx}(v_{x}\bar v +\bar v_{x} v )-2C_{41,x}|v_{x}|^2-C_{41,xxx}|v|^2.
\ea
\ee
Note that it is not hard to verify that $C_{41}$ has the form as in \eqref{a1234}.

Finally, the last two terms $I_{43}$ and $I_{44}$ equal to
\bel{i43}
I_{43}\=a_1(v_{xxxx}\bar v_{x}+\bar v_{xxxx} v_{x})= F_x - a_{1,xxx}|v_{x}|^2+3a_{1,x}|v_{xx}|^2,
\ee
with
$$
F= a_{1}(v_{xxx}\bar v_{x}+\bar v_{xxx} v_{x})-a_{1,x}(v_{xx}\bar v_{x}+\bar v_{xx} v_{x})+a_{1,xx}|v_{x}|^2-a_{1}|v_{xx}|^2,
$$
and
\bel{44}
I_{44}\=a_3(v_{xxxx}\bar v_{xxx}+\bar v_{xxxx} v_{xxx})= \(a_{3}|v_{xxx}|^2\)_x   - a_{3,x}|v_{xxx}|^2.
\ee

By \eqref{11}--\eqref{44}, combining all $``\frac{\partial}{\partial t}$-terms", all $``\frac{\partial}{\partial x}$-terms" and \eqref{pu1} we arrive at the desired inequality \eqref{PWI}.
\section{Global Carleman estimate: Proof of Theorem \ref{Carl}}\label{sec3}

In this section, we obtain a global Carleman estimate inequality for the Schr\"odinger equation \eqref{maineqn} via the poin-wise inequality \eqref{PWI}.
Recalling the definitions of $l$ and $\f$ in \eqref{ellvarphi}, it is easy to check  that
\begin{equation}
\begin{array}{rl}
 |\partial _x^n l|  & \leq C(\psi)\l\mu^n \f, \qq n=1,\cdots, 8,\\
\ns\ns
 |\partial _x^n l_t|& \leq C(\psi)\l\mu^n T\f^2, \q n=1,\cdots, 3,\\
\ns\ns
 |l_t|              & \leq C\l T\f^2,\qq\q |l_{tt}|\leq C\l T^2\f^3.
\end{array}
\end{equation}

We now give the proof of Theorem \ref{Carl}.

{\bf Proof.} The proof is divided into several steps.

{\bf Step $1$}.  Take
$$
\Psi(t,x)=l_x^4.
$$
Recalling the notations in \eqref{a1234}, it is easy to check that the term $\{\cdots\} |v|^2$ in \eqref{PWI} satisfies
\bel{v2}
\{\cdots\} |v|^2=16 l_x^6l_{xx}|v|^2-D_1|v|^2,  \qq |D_1|\geq -C(\psi)\l^6\mu^8\f^6.
\ee
Similarly, we have
\bel{vx2}
\{\cdots\} |v_x|^2=144 l_x^4l_{xx}|v_x|^2-D_2|v_x|^2,  \qq |D_2|\geq -C(\psi)\l^4\mu^6\f^4,
\ee
and
\bel{vxx2}
\{\cdots\} |v_{xx}|^2=24 l_x^2l_{xx}|v_{xx}|^2-D_3|v_{xx}|^2,  \qq |D_3|\geq -C(\psi)\l^2\mu^4\f^2.
\ee
Now we consider those hybrid terms in \eqref{PWI}. It holds
\bel{vvx}
D(v\bar v_x-\bar v v_x) \geq -C(\psi) \l^3\mu^3 T\f^4(|v|^2+|v_x|^2),
\ee
\bel{vvxx}
6il_tl_{xx}(v\bar v_{xx} -\bar vv_{xx} ) \geq -C(\psi) \l^2\mu^2 T\f^3(|v|^2+|v_{xx}|^2),
\ee
\bel{vxvxx}
4il_{tx}(v_{xx}\bar v_x-\bar v_{xx}v_x)\geq -C(\psi) \l\mu T\f^2(|v_x|^2+|v_{xx}|^2),
\ee
\bel{thetav}
\tilde a_0\theta (Pu\bar v+\overline{Pu} v)\geq -C(\psi) \l^4\mu^8 \f^4|v|^2-C(\psi)|\theta Pu|^2,
\ee
and
\bel{thetavxx}
6l_{xx}\theta(Pu\bar v_{xx}+\overline{Pu}v_{xx})\geq -C(\psi) \l^2\mu^4 \f^2|v_{xx}|^2-C(\psi)|\theta Pu|^2.
\ee

Taking \eqref{v2}--\eqref{thetavxx} into \eqref{PWI},  one can find a sufficiently large constant $C(\psi)>0$, only depending on $\psi$, such that
\bel{pointest}
\ba{ll}
&\ds C(\psi)(|\theta P u|^2-A_x-B_t+\l^6\mu^8 \f^6|v|^2+\l^4\mu^6\f^4|v_{x}|^2+\l^2\mu^4\f^2|v_{xx}|^2)\\
\ns\ns\ns
 \geq&\ds 16l_x^6l_{xx}|v|^2+144l_x^4l_{xx}|v_{x}|^2+24l_x^2l_{xx}|v_{xx}|^2+16l_{xx}|v_{xxx}|^2 .
\ea\ee

{\bf Step $2$}. Now we integrate \eqref{pointest} with respect to $t$ and $x$. By the definition of $v=\theta u$ with $\theta(0,x)=\theta(T,x)=0$ and $B$ in \eqref{TermA}, it is obvious that
\bel{BT}
-\int_Q B_t dxdt=0.
\ee
Hence, we have
\bel{temp123}
\ba{ll}
     &\ds  C(\psi)\(\int_Q|\theta P u|^2dxdt-\int_Q A_x dxdt\)\\
           \ns\ns\ns
 \geq&\ds \int_Q(16l_x^6l_{xx}-C(\psi)\l^6\mu^8\f^6)|v|^2dxdt+\int_Q(144l_x^4l_{xx}-C(\psi)\l^4\mu^6\f^4)|v_{x}|^2dxdt\\
            \ns\ns\ns
     &\ds +\int_Q(24l_x^2l_{xx}-C(\psi)\l^2\mu^4\f^2)|v_{xx}|^2dxdt+\int_Q16l_{xx}|v_{xxx}|^2dxdt.
\ea\ee
Since
\bel{lxlxx}
l_x=\l\mu\psi_x\f=\l\mu(x-x_0)\f,\qq l_{xx}=\l\mu(4\mu(x-x_0)^2+2)\f
\ee
by \eqref{x_0} and $\f\leq \frac{T^2}{4}\f^2$, by choosing $\mu\geq \mu_0\geq 1$ and $\l\geq \l_0(\mu)=C(\psi)  (T+T^2)$, it holds that 
\bel{temp11}
 \ds\int_Q(16l_x^6l_{xx}-C(\psi)\l^6\mu^8 \f^6)|v|^2dxdt
\geq  16\int_Q2^8(x-x_0)^8\l^7\mu^8 \f^7 |v|^2dxdt.
 \ee
Similarly,
  \bel{temp12}
  \ds\int_Q(144l_x^4l_{xx}-C(\psi)\l^4\mu^6 \f^4)|v_x|^2dxdt    \geq  144\int_Q2^6(x-x_0)^6\l^5\mu^6 \f^5 |v_x|^2dxdt,
  \ee
  \bel{temp13}
  \ds\int_Q(24l_x^2l_{xx}-C(\psi)\l^2\mu^4 \f^2)|v_{xx}|^2dxdt \geq  24\int_Q2^4(x-x_0)^4\l^3\mu^4 \f^3 |v_{xx}|^2dxdt,
  \ee
and
  \bel{temp14}
  \ds\int_Q16l_{xx}|v_{xxx}|^2dxdt                             \geq  16\int_Q2^2(x-x_0)^2\l  \mu^2 \f |v_{xxx}|^2dxdt.
  \ee
  For the term $A_x$, since $v, v_x, v_t $ and $v_{tx}$ vanish as $x=0,1$ for any $t\in(0,T)$,  we have
$$
\ba{ll}
\ds-\int_Q A_x dxdt=-\int_0^T A(t,1)dt+\int_0^T A(t,0)dt\\
\ns\ns
\ds=\int_0^T \((20 l_x^3+12l_xl_{xx}-10l_{xx})|v_{xx}|^2+4l_x|v_{xxx}|^2+6 l_{xx}(v_{xxx}\bar v_{xx}+\bar v_{xxx} v_{xx})\)\|_{x=0}^{x=1}dt.
\ea
$$
Recalling $l_x$ and $l_{xx}$ in  \eqref{lxlxx},  by taking $\l$ sufficiently large, we have
\bel{AX1}
     \int_0^T A(t,0)dt>0,
\ee
and
\bel{AX2}
     -\int_0^T A(t,1)dt\leq C(\psi)\int_0^T \(\l^3\mu^3 \f^3(t,1)|v_{xx}(t,1)|^2+\l\mu\f(t,1) |v_{xxx}(t,1)|^2\)dt.
\ee
Substituting \eqref{temp11}--\eqref{AX2} into \eqref{temp123}, it holds
\bel{globalest}
\ba{ll}
&\ds \int_Q\l^7\mu^8\f^7|v|^2dxdt+\int_Q\l^5\mu^6\f^5|v_{x}|^2dxdt+\int_Q\l^3\mu^4\f^3|v_{xx}|^2dxdt+\int_Q\l \mu^2\f |v_{xxx}|^2dxdt\\
\ns\ns\ns
 \leq&\ds  C(\psi)\int_Q|\theta P u|^2dxdt+C(\psi)\int_0^T \(\l^3\mu^3 \f^3(t,1)|v_{xx}(t,1)|^2+\l\mu\f(t,1) |v_{xxx}(t,1)|^2\)dt.
\ea\ee
Moreover, since $v=e^{l}u$, we compute
\bel{temp15}
\ba{rl}
  v_{x} &\ds= \theta (u_{x}+l_x u ),\\
 \ns\ns\ns
  v_{xx} &\ds= \theta (u_{xx}+2l_xu_x+(l_x^2+l_{xx})u),\\
\ns\ns\ns
  v_{xxx}  &\ds= \theta(u_{xxx}+3l_x u_{xx}+(3l^2_x +3l_{xx})u_x+(l_x^3+3l_xl_{xx}+l_{xxx})u).
\ea\ee
By Young's inequality, it is not difficult to obtain
\bel{globalest1}
\ba{ll}
&\ds \int_Q\l^7\mu^8\f^7\theta^2|u|^2dxdt+\int_Q\l^5\mu^6\f^5\theta^2|u_{x}|^2dxdt\\
\ns\ns\ns
&\ds +\int_Q\l^3\mu^4\f^3\theta^2|u_{xx}|^2dxdt+\int_Q\l \mu^2\f \theta^2|u_{xxx}|^2dxdt\\
\ns\ns\ns
 \leq&\ds  C(\psi)\int_Q|\theta P u|^2dxdt\\
\ns\ns\ns
     &\ds  +C(\psi)\int_0^T \(\l^3\mu^3 \f^3(t,1)\theta^2(t,1)|u_{xx}(t,1)|^2+\l\mu\f(t,1) \theta^2(t,1)|u_{xxx}(t,1)|^2\)dt,
\ea\ee
which is exactly the statement of Theorem \ref{Carl}.

\vskip 1cm

\section{Boundary observations: Proof of Theorem \ref{invP}}\label{sec4}

In this section, we show the proof of Theorem \ref{invP}, which is a direct application of the Carleman inequality \eqref{Carest}. The standard procedure can be found in \cite{MR1955903, MR2384776}.

{\bf Proof of Th. \ref{invP}.} Pick any $p, q$ as in the statement of the theorem, and introduce the difference $y:=u^p-u^q$ of the corresponding solutions of \eqref{maineqn}.

Then $y$ fulfill the system
  \begin{equation}\label{difeqn}
  \left\{\begin{array}{ll}
  i y_t+y_{xxxx}+q(x)y=f(x)R(t,x), &(t,x)\in Q\\
  y(t,0)=y(t,1)=0, y_x(t,0)=y_x(t,1)=0, & t\in (0,T)\\
y(0,x)=0, & x\in \O.
   \end{array}\right.\end{equation}
   with $f:=q-p$ (real valued) and $R:=u^p$. To complete the proof of Theorem \ref{invP}, we need the following result.

   \begin{proposition}\label{estfx}
Suppose that
\begin{itemize}
\item $R(0, x)\in\dbR$ or $i R(0,x)\in\dbR\;a.e.$ in $\Omega$,
\item $|R(0, x)|\geq r>0\; a.e.$ in $\Omega$,
\item $R\in H^1(0,T;L^\infty(\Omega))$ and
\item $y_{xx}(t,1)\in H^1(0,T)$ and $y_{xxx}(t,1)\in H^1(0,T)$.
\end{itemize}
 Then for any $m\geq 0$ there exists a constant $C>0$ such that for any $q\in L^\infty(\Omega)$ with $\norm{q}_{L^\infty(\Omega)}\leq m$ and for all $f\in L^2(\dbR;\Omega)$, the solution of \eqref{difeqn} satisfies
\begin{equation}\label{inv240}
\norm{f}^2_{L^2(\Omega)}\leq C \(\norm{y_{xx}(\cdot,1)}_{H^1(0,T)}^2+\norm{y_{xxx}(\cdot,1)}_{H^1(0,T)}^2\).
\end{equation}
   \end{proposition}
  {\bf Proof of Proposition \ref{estfx}.} Let $f\in L^2(\dbR; \Omega)$ and $R\in H^1(0,T; L^\infty(\Omega))$ be such that $R(0,x)\in \dbR \; a.e.$ in $\Omega$, and let $y$ be the solution of \eqref{difeqn}. We take the even-conjugate extensions of $y$ and $R$ to the interval $(-T, T)$; i.e., we set $y(t,x)=y(-t,x)$ for $t\in (-T, 0)$ and similarly for $R$. Since $R(0,x)\in \dbR\; a.e.$ in $\Omega$, we have that $R\in H^1(-T,T;L^\infty(\Omega))$, and $y$ satisfies the system \eqref{difeqn} in $(-T, T)\times \Omega$. In the case when $R(0,x)\in i\dbR$, the proof is still valid by take odd-conjugate extensions.

  Changing $t$ into $t+T$, we may assume that $y$ and $R$ are defined on $(0,2T)\times\Omega$, instead of $(-T, T)\times \Omega$.

  Let $z(t,x)=y_t(2T-t, x)$. Then $z$ satisfies the following system:
    \begin{equation}\label{yteqn}
  \left\{\begin{array}{ll}
   z_t+iz_{xxxx}+iq(x)z=if(x)R_t(t,x), &(t,x)\in (0,2T)\times \Omega\\
  z(t,0)=z(t,1)=0, z_x(t,0)=z_x(t,1)=0, & t\in (0,2T)\\
z(T, x)=-if(x)R(T, x), & x\in \O.
   \end{array}\right.\end{equation}
   We shall apply Theorem \ref{Carl}, with $2T$ instead of  $T$. Therefore, here we consider
   $$
   \theta=e^l,\qq\f(t,x)=\ds \frac{e^{3\mu\psi(x)}}{t(2T-t)}\qq\hbox{and}\qq l(t,x)= \l \frac{e^{3\mu\psi(x)}-e^{5\mu\norm{\psi}_\infty}}{t(2T-t)}.
   $$
   To use the Theorem \ref{Carl}, we introduce $v=\theta z$ and $I_1$ is taken as in \eqref{def_i}.

   Now set
   $$
   J=\int_0^T\int_\Omega I_1\theta \bar z dxdt.
   $$
  Then we have
  \begin{equation}
  \begin{array}{lll}
  |J| &\leq&\ds     \(\int_0^T\int_\Omega |I_1|^2dxdt\)^{1/2}
                                     \(\int_0^T\int_\Omega \theta^2|z|^2dxdt\)^{1/2}\\
      \ns\ns
      &\leq&\ds     \lambda^{-7/2}\mu^{-4}
                        \int_0^T\int_\Omega |I_1|^2dxdt
            +           \lambda^{7/2}\mu^{4}
                        \int_0^T\int_\Omega \theta^2|z|^2dxdt\\
      \ns\ns
      &\leq&\ds C   \lambda^{-7/2}\mu^{-4}
                        \(\int_0^T\int_\Omega |I_1|^2dxdt
            +           \lambda^{7}\mu^{8}
                        \int_0^T\int_\Omega \f^7\theta^2|z|^2dxdt\).

\ea\end{equation}
 The last inequality comes from the fact that $\f$ is bounded from below.

 On the other hand, for each $p\in L^\infty(\Omega, \dbR)$,  we define the operator
 $$
 P_p\=\partial_t+i\partial_x^4+ip
 $$
and the space
\bel{zp}
\begin{array}{ll}
&\ds        \mathcal{Z}_p\=    \big\{    z\in L^2(Q);  L_p z\in L^2(Q),  \;z(t,0)=z(t,1)= z_x(t,0)=z_x(t,1)=0, \\
\ns
&   \ds      \qq\qq \hbox{for all}\; \;t\in(0,T),\;     u_{xx}(\cdot, 1)\in L^2(0,T),\; \; u_{xxx}(\cdot, 1)\in L^2(0,T)\big\}.
\ea
\ee
As a direct consequence of Theorem \ref{Carl}, we have the following slightly revised Carleman estimate:
   \bp\label{inv91}
   Given $m\geq0$, there exist $\mu_0\ge1$, $\lambda_0\geq0$ and $C>0$ such that for each $p\in L^\infty(\Omega)$ with $\norm{p}_{L^\infty}\leq m $ it holds
   \be\label{inv237}
\ba{lll}
    \3n &&\3n\ds \int_Q\(\l^7\mu^8\f^7\theta^2|z|^2+|I_1|^2\)dxdt\\
\ns\ns\ns
\3n&\leq &\ds \3n C      \left(   \int_Q\theta^2 |P_p z|^2dxdt
                         +\int_0^T\2n\(\l^3\mu^3\f^3\theta^2|z_{xx}|^2+\l \mu\f \theta^2|z_{xxx}|^2\)(t,1) dt\right)
\ea
\ee
for all $\lambda\geq \lambda_0$, $\mu\geq\mu_0$ and $z\in\mathcal{Z}_p$.
   \ep
{\bf Proof.} The term $|I_1|^2$ can be added by directly taking \eqref{pu1} into account. Moreover, the operator $P$ can be changed to   $P_p$ since $p$ is assumed to be uniformly bounded and the  cost is a slight change of $C$ with respect to the upper bound $m$.

\vskip 3mm

  We now apply the Carleman inequality \eqref{inv237} (with $2T$ instead of $T$) on $z$ and we obtain
    \begin{equation}\label{inv245}
  \begin{array}{lll}
  |J| \3n&\leq&\3n\ds       C   \lambda^{-\frac{7}{2}}\mu^{-4}
                        \(\2n\int_0^{2T}\3n\1n\int_\Omega \theta^2|fR_t|^2dxdt
            +           \int_0^{2T}\3n\(\l^3\mu^3\f^3\theta^2|z_{xx}|^2+\l \mu\f \theta^2|z_{xxx}|^2\)(t,1) dt\)\\
        \ns\ns
      &\leq&\ds C   \lambda^{-\frac{7}{2}}\mu^{-4}
                        \int_\Omega e^{2l(T,x)}|f(x)|^2dx
                             \\
        \ns\ns
      &&\ds +     C   \lambda^{-\frac{1}{2}}\mu^{-1}
                             \int_0^{2T} |z_{xx}(t,1)|^2dt
                   +     C   \lambda^{-\frac{5}{2}}\mu^{-3}
                             \int_0^{2T} |z_{xxx}(t,1)|^2dt.
\ea\end{equation}
  The last inequality holds true due to the fact that $l(T,x)\geq l(t,x)$ for all $(t,x)\in (0, 2T)\times \Omega $, that $\f^3\theta^2 $ and $\f\theta^2$ are bounded from above in $(0, 2T)\times \Omega $ and that $R_t\in L^2(0, 2T; L^\infty(\Omega))$.

   On the other hand,  since $v=\theta z$, we have
\begin{equation*}
J=\int_0^T\int_\Omega I_1 \bar v dxdt\\
=\int_0^T\int_\Omega i v_t\bar v dxdt+ \int_0^T\int_\Omega \((\Psi+\frac{1}{2} a_{2,xx}) |v|^2-a_2 |v_x|^2 +| v_{xx}|^2\)dxdt,
\end{equation*}
hence,
\begin{equation*}
\Im (J)=\frac{1}{2}\int_\Omega |v(T, x)|^2 dx=\frac{1}{2}\int_\Omega  e^{2l(T, x)} |f(x)|^2 |R(T, x)|^2 dx.
\end{equation*}
Using the hypothesis on $R(T, x)$, it  follows that
\begin{equation}\label{inv243}
\Im (J)\geq \frac{r^2}{2}\int_\Omega  e^{2l(T, x)} |f(x)|^2  dx.
\end{equation}

  Combining \eqref{inv245} and \eqref{inv243}, we have that
  \bel{inv246}
 \int_\Omega  e^{2l(T, x)} |f(x)|^2  dx\leq C  \(\int_0^{2T} |z_{xx}(t,1)|^2
                   +     |z_{xxx}(t,1)|^2\)dt
  \ee
  for $\lambda $ and $\mu$ large enough. Then \eqref{inv240} follows from \eqref{inv246} since
  $$
  e^{2l(T,x)}\geq e^{2M}>0, \qq \hbox{with}\qq \ds M=\frac{\lambda}{T^2}(1-e^{5\mu\norm{\psi}_{\infty}}).
  $$
  This completes the proof of Proposition \ref{estfx} and of Theorem \ref{invP}.

  \section{Further comments and open problems}

  \begin{enumerate}

\item  There is another formulation for stationary inverse problems known as the Dirichlet-to-Neumann map. For instance,   Bukhgeim and Uhlmann (\cite{MR1900557}) show that the potential can be uniquely determined by  the boundary data for  $(\Delta+q)u=0$.  It would be interesting to find out what happens for the fourth order Schr\"odinger operator. However, the relationship between the two problems is not really clear.

\item   In this paper we derive a boundary Carleman estimate for the fourth order Schr\"odinger operator.  It is well known that  based on \eqref{Carest}, we can derive the observability inequality and, consequently,  prove the controllability property of  the controlled system with two boundary controls.   As a direct consequence of this methodology, it is very likely to expect that the controllability property holds for the fourth order Schr\"odinger equation with nontrivial potential $q$.  Such result is much more general than the existing one in \cite{MR2928966},  which is for trivial potential $q$,  even though only one boundary control is needed.  It would be interesting to know whether two controls on the boundary are necessary with the nontrivial potential $q$.

\item  It is well known that  the Carleman estimate is a useful tool to analyze inverse problems. In  fact, it has been studied for second order Schr\"odinger operator not only in bounded domain, but also  in an unbounded strip (\cite{MR2394028}) or on a tree (\cite{MR2864510}). One could expect similar results in different domains.   Meanwhile, it is still a challenging problem whether one can construct Carleman inequalities  for fourth order equations on higher dimensions.

\item Note that there are fruitful literatures considering the numerical approximation results for the second order Schr\"odinger equations. Similar to the discrete Carleman estimate constructed by parabolic equation (see \cite{MR2745778}),  it would be interesting  to find out the discrete analogue of \eqref{Carest} for space semi-discretized Schr\"odinger equation as the first step to solve  discrete problems.

\end{enumerate}

\end{document}